\documentclass{amsart}
\usepackage{amssymb, amsthm, amscd, amsmath}

\newtheorem{theorem}{Theorem}[section]
\newtheorem{lemma}[theorem]{Lemma}
\newtheorem{prop}[theorem]{Proposition}

\theoremstyle{definition}
\newtheorem{defn}[theorem]{Definition}
\newtheorem{remark}[theorem]{Remark}
\newtheorem{question}[theorem]{Question}

\def\Ker{\mathop{\mathrm{Ker}}}
\def\Im{\mathop{\mathrm{Im}}}
\def\ann{\mathop{\mathrm{ann}\,}}
\def\Ass{\mathop{\mathrm{Ass}}}
\def\Pdot{\mathop{\mathrm{P_{\bullet}}}}
\def\Min{\mathop{\mathrm{Min}}}
\def\Spec{\mathop{\mathrm{Spec}}}
\def\Supp{\mathop{\mathrm{Supp}}}
\def\m{\mathfrak{m}}
\def\p{\mathfrak{p}}
\def\P{\mathfrak{P}}
\def\q{\mathfrak{q}}
\def\a{\mathfrak{a}}
\def\N{\mathbb{N}}
\def\Z{\mathbb{Z}}

\begin{document}

\title[Associated primes]{Associated primes of local cohomology modules and
of Frobenius powers}

\author{Anurag K. Singh \and Irena Swanson}

\address{
School of Mathematics \\
686 Cherry St. \\
Georgia Institute of Technology \\
Atlanta, GA 30332-0160, USA.
E-mail: {\tt singh@math.gatech.edu}
\phantom{.} \newline
\phantom{} \quad Department of Mathematical Sciences \\
New Mexico State University \\
Las Cruces, NM~88003-8001, USA.
E-mail: {\tt iswanson@nmsu.edu}}

\thanks {Both authors were supported in part by grants from the National
Science Foundation.}


\dedicatory{Dedicated to Professor Melvin Hochster on the occasion of his
sixtieth birthday}

\date{\today}

\begin{abstract}
We construct normal hypersurfaces whose local cohomology modules have
infinitely many associated primes. These include unique factorization domains
of characteristic zero with rational singularities, as well as F-regular unique
factorization domains of positive characteristic. As a consequence, we answer a
question on the associated primes of Frobenius powers of ideals, which arose
from the localization problem in tight closure theory.
\end{abstract}

\maketitle

\section{Introduction}

Let $R$ be a commutative Noetherian ring and ${\mathfrak a} \subset R$ an 
ideal. In \cite{Hu1} Huneke asked whether the number of associated prime ideals
of a local cohomology module $H_{\mathfrak a}^n(R)$ is always finite. In
\cite{Si} the first author constructed an example of a hypersurface
$$
R = {\Z}[u,v,w,x,y,z]/(ux+vy+wz)
$$
for which the local cohomology module $H_{(x,y,z)}^3(R)$ has a $p$-torsion
element for every prime integer $p$, and consequently has infinitely many 
associated prime ideals. However this example does not address Huneke's 
question for rings containing a field, nor does it yield an example over a 
local ring.
More recently Katzman constructed the following example in \cite{Ka2}:
let $K$ be an arbitrary field, and consider the hypersurface
$$
S = K[s,t,u,v,x,y]/\big(su^2x^2-(s+t)uxvy+tv^2y^2\big).
$$
Katzman showed that the local cohomology module $H^2_{(x,y)}(S)$ has infinitely
many associated prime ideals. Since the defining equation of this hypersurface
factors, the ring in Katzman's example is not an integral domain. In this paper
we generalize Katzman's construction and obtain families of examples which
include examples over normal domains, and even over hypersurfaces with rational
singularities:

\begin{theorem}\label{lcintro}
Let $K$ be an arbitrary field. Then there exists a standard graded hypersurface
$R$ with $[R]_0=K$, which is a unique factorization domain, and contains an
ideal $\a$, such that a local cohomology module $H^n_{\a}(R)$ has infinitely
many associated prime ideals.

If $K$ has characteristic zero, there exist such examples where, furthermore,
$R$ has rational singularities. If $K$ has positive characteristic, we may
choose $R$ to be F-regular. If $\m$ denotes the homogeneous maximal ideal of
$R$, then $H^n_{\a}\big(R_\m\big)$ has infinitely many associated prime ideals
as well.
\end{theorem}

There are affirmative answers to Huneke's question if the ring $R$ is regular
but, as our theorem indicates, the hypothesis of regularity cannot be weakened
substantially. The first results were obtained by Huneke and Sharp who proved
that if $R$ is a regular ring containing a field of prime characteristic, then
the set of associated prime ideals of $H_{\mathfrak a}^n(R)$ is finite,
\cite[Corollary 2.3]{HS}. Lyubeznik established that $H_{\mathfrak a}^n(R)$ has
finitely many associated prime ideals if $R$ is a regular local ring containing
a field of characteristic zero, or an unramified regular local ring of mixed
characteristic, see \cite[Corollary 3.6\,(c)]{Ly1} and \cite[Theorem 1]{Ly2}
respectively. Marley proved that if $R$ is a local ring, then for any finitely
generated $R$-module $M$ of dimension at most three, any local cohomology module
$H_{\mathfrak a}^n(M)$ has finitely many associated primes, \cite[Corollary
2.7]{Ma}. If $i$ is the smallest integer for which $H_\a^i(M)$ is not a
finitely generated $R$-module, then the set $\Ass H_\a^i(M)$ is finite, as
proved in \cite{BF} and \cite{KS}. For some of the other work on this question,
we refer the reader to the papers \cite{BKS, BRS, He, Ly3, MV} and \cite{TZ}.

In \S\ref{froblc} we establish a relationship between the associated primes of
Frobenius powers of an ideal and the associated primes of a local cohomology
module over an auxiliary ring. Recall that for an ideal $\a$ in a ring $R$ of
prime characteristic $p>0$, the {\it Frobenius powers}\/ of $\a$ are the ideals
$\a^{[p^e]} = \left(x^{p^e} \ | \ x \in \a\right)$ where $e \in \N$. The
finiteness of the associated primes of the ideals $\a^{[p^e]}$ is related to
the localization problem in tight closure theory, discussed in \S\ref{tc} of
this paper. In \cite{Ka1} Katzman constructed the first example where the set
$\bigcup_e\Ass R/\a^{[p^e]}$ is infinite. The question however remained
whether the set $\bigcup_e\Ass R/\big(\a^{[p^e]}\big)^*$ is finite, or if it
has finitely many maximal elements---this has strong implications for the
localization problem, see \cite{AHH, Ho, Ka1, SN} or \cite[\S12]{Hu2}. As an
application of our results on local cohomology, we settle this question in
\S\ref{tc}, with the following theorem:

\begin{theorem}
There exists an F-regular unique factorization domain
$R$ of characteristic $p>0$, with an ideal $\a$,
for which the set
$$
\bigcup_{e \in \N}\Ass \frac{R}{\a^{[p^e]}}
= \bigcup_{e \in \N}\Ass \frac{R}{\big(\a^{[p^e]}\big)^*}
$$
has infinitely many maximal elements.
\end{theorem}

\section{General constructions}\label{froblc}

Let $\a=(x_1,\dots,x_n)$ be an ideal of a ring $R$. For an integer $r \ge 0$,
the local cohomology module $H^r_{\a}(R)$ may be computed as the $r$th
cohomology module of the extended \v Cech complex
$$
0 \longrightarrow R \longrightarrow \bigoplus_{i=1}^n R_{x_i} \longrightarrow 
\bigoplus_{i<j} R_{x_i x_j} \longrightarrow \cdots \longrightarrow
R_{x_1 \cdots x_n} \longrightarrow 0.
$$
For positive integers $m_i$ and an element $f \in R$, we will use 
$[f+(x_1^{m_1},\dots,x_n^{m_n})]$ to denote the cohomology class
$$
\left[\frac{f}{x_1^{m_1} \cdots x_n^{m_n}} \right] \ \in \ H^n_{\a}(R) \ 
= \ \frac{R_{x_1 \cdots x_n}}{\sum R_{x_1 \cdots \hat{x_i} \cdots x_n}}.
$$
It is easily seen that $[f+(x_1^{m_1},\dots,x_n^{m_n})] \in H^n_{\a}(R)$ is
zero if and only if there exist integers $k_i \ge 0$ such that
$$
fx_1^{k_1}\cdots x_n^{k_n} \ \in\ \big(x_1^{m_1+k_1},\dots,x_n^{m_n+k_n}\big)R. 
$$
Consequently $H^n_{\a}(R)$ may also be computed as the direct limit
$$
H^n_{\a}(R) \cong \varinjlim_{m \in \N} R/(x_1^m,\dots,x_n^m)R,
$$
where the maps in the direct system are induced by multiplication by the
element $x_1 \cdots x_n$. We may regard an element $[f+(x_1^m,\dots,x_n^m)] \in
H^n_{\a}(R)$ as the class of $f+(x_1^m,\dots,x_n^m)R$ in this direct limit.

We next record two results which illustrate the relationship between associated
primes of local cohomology modules and associated primes of generalized
Frobenius powers of ideals. 

\begin{prop}\label{contain}
Let $R$ be a Noetherian ring, and $\{M_i\}_{i\in I}$ be a direct system of
$R$-modules. Then
$$
\Ass \left( \varinjlim M_i \right) \subseteq \bigcup_{i\in I}\Ass M_i.
$$
In particular, if $\a=(x_1,\dots,x_n)$ is an ideal of $R$, then for any infinite
set $\mathbb S$ of positive integers,
$$
\Ass H^n_{\a}(R) \subseteq \bigcup_{m \in \mathbb S} 
\Ass \frac{R}{(x_1^m, \ldots, x_n^m)}.
$$
\end{prop}

\begin{proof}
Let $\p = \ann m$ for some element $m \in \varinjlim M_i$. If $z \in \p$
then $zm = 0$, and so there exists $i \in I$ such that $m$ is the image of
$m_i \in M_i$ and $zm_i = 0$. Since $\p$ is finitely generated, there exists $j
\ge i$ such that $m_i \mapsto m_j \in M_j$ and $\p m_j = 0$. Consequently
$$
\p \subseteq 0 :_R m_j \subseteq 0 :_R m = \p,
$$
i.e, $\p = {\ann} m_j \in\Ass M_j$. 
\end{proof}

It immediately follows that whenever $H^n_{\a}(R)$ has infinitely many
associated prime ideals, the set $\bigcup_m\Ass R/(x_1^m, \ldots, x_n^m)$ is
infinite as well. The converse is false, as we shall see in
Remark~\ref{examplelcfrob}.

\begin{prop} \label{multigraded}
Let $A$ be an $\N$-graded ring which is generated, as an $A_0$-algebra, by
elements $t_1,\dots,t_n$ of degree $1$ which are nonzerodivisors in $A$. Let
$R$ be the extension ring
$$
R = A[u_1,\dots,u_n,x_1,\dots,x_n]/(u_1x_1-t_1,\dots,u_nx_n-t_n).
$$
Let $m_1,\dots,m_n$ be positive integers, and $f \in A$ a homogeneous element.
Then, for arbitrary integers $k_i \ge 0$,
$$
\big(t_1^{m_1},\dots,t_n^{m_n}\big) A :_{A_0} f \ = \
\big(x_1^{m_1+k_1},\dots,x_n^{m_n+k_n}\big)R :_{A_0} fx_1^{k_1}\cdots x_n^{k_n}.
$$
Consequently, if we consider the element $\eta=[f+(x_1^{m_1},\dots,x_n^{m_n})]$
of the local cohomology module $H^n_{(x_1, \dots, x_n)}(R)$, then
$$
\big(t_1^{m_1},\dots,t_n^{m_n}\big)A :_{A_0} f \ = \ {\ann}_{A_0}\eta.
$$
\end{prop}

\begin{proof}
The inclusion $\subseteq$ is easily verified. For the other inclusion, let 
$e_i \in {\Z}^{n+1}$ be the unit vector with $1$ as its $i$th entry, and
consider the ${\Z}^{n+1}$-grading on $R$ where $\deg x_i=e_i$ and
$\deg u_i=e_{n+1}-e_i$ for all $1\le i\le n$. If $f\in A_r$ then, as an element
of $R$, the degree of $f$ is $re_{n+1}$. The subring $A$ is a direct summand of
$R$ since
$$
A_j=R_{(0,\dots,0,j)} \quad \text{for} \quad j \ge 0, \quad \text{and} \quad
A=\oplus_{j \ge 0} R_{(0,\dots,0,j)}.
$$
Now if $h\in A_0$ is an element such that
$hfx_1^{k_1}\cdots x_n^{k_n} \in (x_1^{m_1+k_1},\dots,x_n^{m_n+k_n})R$, then
there exist homogeneous elements $c_1,\dots,c_n \in R$ such that
$$
hfx_1^{k_1} \cdots x_n^{k_n} \ = \ c_1x_1^{m_1+k_1}+\dots+c_nx_n^{m_n+k_n}.
$$
Comparing degrees, we must have $\deg c_1=(-m_1,k_2,\dots,k_n,r)$, and so $c_1$
is an $A_0$-linear combination of monomials $\mu$ of the form
$$
\mu = u_1^{l_1+m_1} \ u_2^{l_2} \ \cdots \ u_n^{l_n} \ 
x_1^{l_1} \ x_2^{l_2+k_2} \ \cdots \ x_n^{l_n+k_n},
$$
where $l_i \ge 0$, and $m_1+l_1+\dots+l_n=r$. Consequently
\begin{align*}
\mu x_1^{m_1+k_1} \ &= \ (u_1x_1)^{l_1+m_1} \ (u_2x_2)^{l_2} \ \cdots \ 
(u_nx_n)^{l_n} \ x_1^{k_1} \ \cdots \ x_n^{k_n} \\
&= \ t_1^{l_1+m_1} \ t_2^{l_2} \ \cdots \ t_n^{l_n} \ x_1^{k_1} \ \cdots \
x_n^{k_n},
\end{align*}
and so $c_1x_1^{m_1+k} \in \big(x_1^{k_1} \cdots x_n^{k_n}t_1^{m_1}\big)R$.
Similar computations for $c_2,\dots,c_n$ show that
$$
hfx_1^{k_1} \cdots x_n^{k_n} \ \in \ x_1^{k_1} \cdots x_n^{k_n}
\big(t_1^{m_1},\dots,t_n^{m_n}\big)R.
$$
Multiplying by $u_1^{k_1} \cdots u_n^{k_n}$ and using that $A$ is a direct
summand of $R$, we get 
\begin{align*}
hf t_1^{k_1} \cdots t_n^{k_n} & \ \in \
t_1^{k_1} \cdots t_n^{k_n} \big(t_1^{m_1},\dots,t_n^{m_n}\big)R \cap A \\
&= \ t_1^{k_1} \cdots t_n^{k_n} \big(t_1^{m_1},\dots,t_n^{m_n}\big)A.
\end{align*}
Since the elements $t_i \in A$ are nonzerodivisors, the required result follows.
\end{proof}

We next record two results which will be used in the proof of
Theorem~\ref{thmlc}.

\begin{lemma}\label{minprimes}
Let $M$ be a square matrix with entries in a ring $R$. Then the minimal primes
of the ideal $(\det M)R$ are precisely the minimal primes of the cokernel of
the matrix $M$.
\end{lemma}

\begin{proof}
Let $C$ denote the cokernel of $M$, i.e., we have an exact sequence
$$
R^n \overset{M}\longrightarrow R^n \longrightarrow C \longrightarrow 0.
$$
For a prime ideal $\p \in \Spec R$, note that $C_{\p}=0$ if and only if
$R_\p^n \overset{M}\longrightarrow R_\p^n$ is surjective or, equivalently, is an
isomorphism. This occurs if and only if $\det M$ is a unit in $R_{\p}$, and so
we have
$$
V(\det M) = \Supp C.
$$
\end{proof}

\begin{lemma}\label{subring}
Let $R$ be an $\N$-graded ring, and $M$ be a $\Z$-graded $R$-module. For every
integer $r$ and prime ideal $\p \in {\Ass}_{R_0}M_r$, there exists a
homogeneous prime ideal $\P \in {\Ass}_R M$ such that $\P \cap R_0=\p$.
Consequently, if the set ${\Ass}_{R_0} M$ is infinite, then so is the set
${\Ass}_R M$.
\end{lemma}

\begin{proof}
Let $\p={\ann}_{R_0}m$ for some element $m \in M_r$. There is no loss of
generality in replacing $M$ by the cyclic module $R/\a \cong mR$, in which case
$\p=\a \cap R_0$. The isomorphism
$$
R/(\a+R_+) \cong R_0 / \p
$$
shows that $\a+R_+$ is a prime ideal of $R$. Let $\P$ be a minimal prime of
$\a$ which is contained in $\a+R_+$. Then $\P \in{\Min}_R R/\a \subseteq 
{\Ass}_R R/\a$, and $\P \cap R_0 = \p$ since $(\a+R_+) \cap R_0 = \p$.
\end{proof}

\begin{defn}\label{multidiag}
Let $d$ be a positive even integer, and $r_0,\dots,r_d$ be elements of a ring
$A_0$. The $n$th {\it multidiagonal matrix with respect to $r_0,\dots,r_d$}
will refer to the $n \times n$ matrix
$$
M_n=\left[\begin{matrix}
r_{\frac{d}{2}} & \dots & r_0 \cr
\vdots & \ddots & & \ddots \cr
r_d & & \ddots & & \ddots \cr
& \ddots & & \ddots & & r_0 \cr
& & \ddots & & \ddots & \vdots \cr
& & & r_d & \dots & r_{\frac{d}{2}} \cr
\end{matrix}\right],
$$
where the elements $r_0,\dots,r_d$ occur along the $d+1$ central diagonals, and
all the other entries are zero. (These multidiagonal matrices are special cases
of T\"oplitz matrices.)
\end{defn}

\begin{theorem}\label{thmlc}
Let $d$ be an even positive integer, $r_0,\dots,r_d$ elements of a domain $A_0$,
$a\ge 0$ an integer, and $M_n$ the $n$th multidiagonal matrix with respect to
$r_0,\dots,r_d$. Let $u,v,x,y$ be variables over $A_0$, and
$\mathbb S \subseteq \N$ a subset such that
$$
\bigcup_{n \in \mathbb S}{\Min}\left(\det M_{n-a-d/2}\right)
$$
is an infinite set. If
$$
A = A_0[x,y]/(xy)^a\left(r_0 x^d+r_1x^{d-1}y+\dots+r_d y^d\right),
$$
then $\bigcup_{n \in \mathbb S}\Ass A/(x^n,y^n)$, is an infinite set.

Furthermore, if $r_0$ and $r_d$ are nonzero elements of $A_0$, then for
$$
R=A_0[u,v,x,y]/\left(r_0(ux)^d+r_1(ux)^{d-1}(vy)+\dots+r_d(vy)^d\right),
$$
the local cohomology module $H^2_{(x,y)}(R)$ has infinitely many associated
primes.

If $(A_0,\m)$ is a local domain or if $(A_0,\m)$ is a graded domain and $\det
M_n$ is a homogeneous element for all $n\ge 0$, then these issues are preserved
under localizations of $A$ and $R$ at the respective maximal ideals
$(\m+(x,y))A$ and $(\m+(u,v,x,y))R$.
\end{theorem}

\begin{proof} 
Consider the $A_0$-module $[A/(x^n,y^n)]_{n-1+a+d/2}$ for $n>a+d$. A generating
set for this module is given by the $n-a-d/2$ monomials
$$
x^{a+d/2}y^{n-1}, \ x^{a+d/2+1}y^{n-2}, \ \dots, \ x^{n-1}y^{a+d/2}.
$$
There are $n-a-d/2$ relations amongst these monomials, arising from the
equations
$$
(xy)^a\big(r_0x^d+r_1x^{d-1}y+\dots+r_d y^d\big)x^iy^{n-1-a-d/2-i} = 0,
$$
where $0 \le i \le n-1-a-d/2$. Using this, it is easily checked that the
presentation matrix for $[A/(x^n,y^n)]_{n-1+a+d/2}$ is precisely the
multidiagonal matrix $M_{n-a-d/2}$. By Lemma~\ref{minprimes}, whenever 
$\det M_{n-a-d/2}$ is nonzero, its minimal primes are the minimal primes of
$[A/(x^n,y^n)]_{n-1+a+d/2}$, and so
$$
\bigcup_{n\in\mathbb S}{\Ass}_{A_0} \left[\frac{A}{(x^n,y^n)}\right]_{n-1+a+d/2}
$$
is an infinite set. Using Lemma~\ref{subring}, the set
$\bigcup_{n\in\mathbb S}\Ass A/(x^n,y^n)$ is infinite as well. 

Note that $xy$ is a nonzerodivisor in $A_0[x,y]/(r_0 x^d + r_1x^{d-1}y +
\dots+r_d y^d)$ whenever $r_0$ and $r_d$ are nonzero elements of $A_0$. The set
${\Ass}_{A_0} H^2_{(x,y)}(R)$ is infinite by Proposition~\ref{multigraded}.
Since $A_0 = R_0$, Lemma~\ref{subring} implies that the set ${\Ass}_R
H^2_{(x,y)}(R)$ is infinite.
\end{proof}

\begin{remark}\label{moty}
We demonstrate how Katzman's examples from \cite{Ka1} and \cite{Ka2} follow
from Theorem~\ref{thmlc}. Let $K$ be an arbitrary field, and consider the
polynomial ring $A_0=K[t]$. Let $M_n$ be the $n$th multidiagonal matrix with 
respect to the elements $r_0=1$, $r_1=-(1+t)$, and $r_2=t$. An inductive
argument shows that
$$
\det M_n = (-1)^n(1+t+t^2+\cdots+t^n) = (-1)^n\frac{t^{n+1}-1}{t-1}
\quad \text{for all} \quad n \ge 1.
$$
It is easily verified that $\bigcup_{n \in \N}{\Min}(\det M_n)$ is an infinite
set and, if $K$ has characteristic $p>0$,
that the set $\bigcup_{e \in \N}{\Min}(\det M_{p^e-2})$ is also infinite.
Theorem~\ref{thmlc}
now gives us the main results of \cite{Ka2}: the local cohomology module
$H^2_{(x,y)}(R)$ has infinitely many associated primes where
$$
R = K[t,u,v,x,y]/\big(u^2x^2-(1+t)uxvy+tv^2y^2\big).
$$
Similarly, graded or local examples may be obtained using
$$
S = K[s,t,u,v,x,y]/\big(su^2x^2-(s+t)uxvy+tv^2y^2\big),
$$
in which case $H^2_{(x,y)}(S)$ and $H^2_{(x,y)}(S_\m)$ have infinitely many
associated primes.

If $K$ has characteristic $p>0$, consider the hypersurface
$$
A = K[t,x,y]/\big(xy(x^2-(1+t)xy+ty^2) \big),
$$
where $a=1$ in the notation of Theorem~\ref{thmlc}. The theorem now implies
that the Frobenius powers of the ideal $(x,y)A$ have infinitely many associated
primes, as first proved by Katzman in \cite{Ka1}.
\end{remark}

\section{Tridiagonal matrices}\label{d=2}

The results of the previous section demonstrate how multidiagonal matrices give
rise to associated primes of local cohomology modules and of Frobenius powers of
ideals. One of the goals of this paper is to construct an integral domain $A$ of
characteristic $p>0$, with an ideal $\a$, such that the set
$\bigcup_e\Ass A/\a^{[p^e]}$ is infinite. To obtain such examples directly from
Theorem~\ref{thmlc} we need the set
$\bigcup_e{\Min}\left(\det M_{p^e-d/2}\right)$ to be infinite, since the domain
hypothesis forces $a=0$ in the notation of the theorem. In \S\ref{d=4} we show
that $\bigcup_e{\Min}\left(\det M_{p^e-d/2}\right)$ can indeed be infinite when
$d=4$.

In Proposition~\ref{threediag} we prove that
$\bigcup_e{\Min}\left(\det M_{p^e-d/2}\right)$ is finite whenever $d=2$, see
also \cite[Lemma 10]{Ka1}. Nevertheless, the main results of our paper rely 
heavily on an analysis of multidiagonal matrices with $d=2$, which we undertake
next.

In the notation of Definition~\ref{multidiag}, multidiagonal matrices with $d=2$
have the form
$$
M_n=\left[\begin{matrix}
r_1 & r_0 \cr
r_2 & r_1 & r_0 \cr
& \ddots & \ddots & \ddots \cr
& & r_2 & r_1 & r_0 \cr
& & & r_2 & r_1 \cr
\end{matrix}\right].
$$
It is convenient to define $\det M_0=1$, and it is easily seen that
$$
\det M_{n+2}=r_1 \det M_{n+1}-r_0r_2 \det M_n \quad \text{for all} \quad n \ge 0.
$$
While we will not be using it here, we mention that
$$
\det M_n=\sum_{i=0}^{\lfloor n/2\rfloor}(-1)^i \binom{n-i}{i}
r_1^{n-2i}(r_0r_2)^i.
$$
Consider the generating function for $\det M_n$, 
$$
G(x) = \sum_{n \ge 0} \left(\det M_n\right) x^n.
$$
By the recursion formula,
$$
\sum_{n \ge 0} \left(\det M_{n+2}\right)x^{n+2} 
= r_1 \sum_{n \ge 0} \left(\det M_{n+1}\right) x^{n+2}
- r_0r_2 \sum_{n \ge 0} \left(\det M_n\right) x^{n+2}
$$
and substituting $G(x)$ and solving, we get
$$
G(x) = \sum_{n \ge 0} \left(\det M_n\right)x^n = \frac{1}{1-r_1 x+r_0r_2 x^2}.
$$

\begin{prop} \label{threediag}
Let $r_0,r_1,r_2$ be elements of a ring $R$ of prime characteristic $p>0$. For
each $n \in \N$, let $M_n$ be the $n$th multidiagonal matrix with respect to
$r_0,r_1,r_2$. Then, for any integer $e \ge 1$,
$$
\det M_{p^e-1} = \left(\det M_{p-1}\right)^{1+p+\dots+p^{e-1}}.
$$
Consequently, the set $\bigcup_e{\Min}\left(\det M_{p^e-1}\right)$ is finite.
\end{prop}

\begin{proof}
Let $1-r_1x+r_0r_2x^2 = (1-\alpha x)(1-\beta x)$ for some elements $\alpha$ and
$\beta$ in a suitable extension of $R$. The generating function $G(x)$ can be
written as
$$
G(x) = \sum_{n \ge 0} \left(\det M_n\right) x^n
= \frac{1}{(1-\alpha x)(1-\beta x)} = \sum_{i,j \ge 0} \alpha^i \beta^j x^{i+j},
$$
and consequently
$$
\det M_{p-1} = \sum_{i=0}^{p-1} \alpha^i \beta^{p-1-i} \qquad \text{and} \text
\qquad \det M_{p^e-1} = \sum_{i=0}^{p^e-1} \alpha^i \beta^{p^e-1-i}.
$$
Using this,
\begin{multline*}
(\det M_{p-1})^{1+p+\cdots+p^{e-1}}
= \prod_{j=0}^{e-1} \left(\sum_{i=0}^{p-1} \alpha^i \beta^{p-1-i}\right)^{p^j}
= \prod_{j=0}^{e-1} \left(\sum_{i=0}^{p-1} \alpha^{ip^j} \beta^{(p-1-i)p^j} 
\right) \\
= \sum_{k=0}^{p^e-1} \alpha^k \beta^{p^e-1-k} = \det M_{p^e-1}.
\end{multline*}
\end{proof}

We next consider a special family of tridiagonal matrices: let $K[s,t]$ be a
polynomial ring over a field $K$, and consider the $n \times n$ multidiagonal
matrices
$$
M_n=\left[\begin{matrix}
t & s \cr
s & t & s \cr
& \ddots & \ddots & \ddots \cr
& & s & t & s \cr
& & & s & t \cr
\end{matrix}\right].
$$
In the notation of Definition~\ref{multidiag}, we have $d=2$, $r_1=t$, and 
$r_0=r_2=s$. Setting $Q_n(s,t) = \det M_n$, we have
$$
Q_0=1, \quad Q_1=t, \quad \text{and} \quad Q_{n+2}=tQ_{n+1}-s^2Q_n \quad
\text{for all} \quad n \ge 0.
$$
Note that the polynomials $Q_n(s,t)$ are relatively prime to $s$. Using the 
specialization $P_n(t)=Q_n(1,t)$, we get polynomials $P_n(t) \in K[t]$ 
satisfying the recursion
$$
P_0(t)=1, \quad P_1(t)=t, \quad \text{and} \quad P_{n+2}(t)=tP_{n+1}(t)-P_n(t)
\quad \text{for all} \quad n \ge 0.
$$
Each $P_n(t)$ is a monic polynomial of degree $n$, and in Lemma~\ref{inf} we 
establish that the number of distinct irreducible factors of the polynomials 
$\{P_n(t)\}_{n \in \N}$ is infinite. As $Q_n(s,t)=s^nP_n(t/s)$ for all $n \ge
0$, this also establishes that the number of distinct irreducible factors of
the polynomials $\{Q_n(s,t)\}_{n \in \N}$ is infinite.

\begin{lemma}\label{roots}
Let $K$ be an algebraically closed field and consider the polynomials 
$P_n(t)=\det M_n \in K[t]$ for $n \ge 1$ as above.
\begin{enumerate}
\item If $\xi$ is a nonzero element of $K$ with $\xi \neq \pm 1$, then 
$P_n(\xi+\xi^{-1})=0$ if and only if $\xi^{2n+2}=1$.

\item The number of distinct roots of $P_n$ which are different from $0$ and 
$\pm 1$ is half of the number of distinct $(2n+2)$th roots of unity different
from $\pm 1$.

\item If $2n+2$ is invertible in $K$, then $P_n(t)$ has $n$ distinct roots of
the form $\xi+\xi^{-1}$ where $\xi^{2n+2}=1$ and $\xi \neq \pm 1$.

\item The elements $2$ or $-2$ are roots of $P_n(t)$ if and only if the
characteristic of $K$ is a positive prime $p$ which divides $n+1$.

\item If the characteristic of $K$ is an odd prime $p$, then $P_{q-2}(t)$ has 
$q-2$ distinct roots for all $q=p^e$. If $p=2$, then $P_{q-2}(t)$ has $q/2-1$
distinct roots.
\end{enumerate}
\end{lemma}

\begin{proof}
(1) Consider the generating function of the polynomials $P_n(t)$, 
$$
G(t,x) = \sum_{n \ge 0} P_n(t)x^n = \frac{1}{1-xt+x^2} \in K[t][[x]].
$$
If $\xi\neq 0$ and $\xi \neq \pm 1$, then
\begin{align*}
\sum_{n \ge 0} P_n(\xi+\xi^{-1})x^n & = \frac{1}{1-x(\xi+\xi^{-1})+x^2} 
= \frac{1}{(\xi^{-1}-x)(\xi-x)} \\
& = \frac{1}{(\xi-\xi^{-1})(\xi^{-1}-x)} - \frac{1}{(\xi-\xi^{-1})(\xi-x)} \\
& = \frac{\xi}{\xi-\xi^{-1}} \sum_{n \ge 0} (\xi x)^n
- \frac{\xi^{-1}}{\xi-\xi^{-1}} \sum_{n \ge 0} (\xi^{-1} x)^n \in K[[x]].
\end{align*}
Equating the coefficients of $x^n$, we have
$$
P_n(\xi+\xi^{-1}) = \frac{\xi^{n+1}-\xi^{-(n+1)}}{\xi-\xi^{-1}}
= \frac{\xi^{2n+2}-1}{\xi^n(\xi^2-1)},
$$
and the assertion follows.

(2) We observe that 
$$
\xi+\frac{1}{\xi} - \left(\eta+\frac{1}{\eta}\right) 
= \xi-\eta - \frac{\xi - \eta}{\xi \eta}
= (\xi-\eta)\left(1 - \frac{1}{\xi \eta}\right),
$$
and so $\xi+\xi^{-1} = \eta+\eta^{-1}$ if and only if $\xi$ equals $\eta$ or
$\eta^{-1}$. 

(3) Since $2n+2$ is invertible in $K$, the polynomial $X^{2n+2}-1=0$ has $2n$
distinct roots $\xi$ with $\xi \neq \pm 1$. These give the $n$ distinct roots
$\xi+\xi^{-1}$ of the degree $n$ polynomial $P_{n}(t)$.

(4) Using the generating function above, 
$$
G(2,x) = \frac{1}{1-2x+x^2} = (1-x)^{-2} = 1+2x+3x^2+\cdots,
$$
and so $P_n(\pm 2)=0$ if and only if $n+1=0$ in $K$.

(5) The case when $p$ is odd follows immediately from (2). If $p=2$, the
equation $X^{2q-2}-1=(X^{q-1}-1)^2=0$ has $q-2$ distinct roots $\xi$ with $\xi
\neq 1$, which ensures that $P_{q-2}(t)$ has at least $q/2-1$ distinct roots.
It follows from (4) that $0$ is not a root of $P_{q-2}(t)$, so these must be
all the roots.
\end{proof}

\begin{lemma}\label{inf}
Let $K$ be an arbitrary field. Then the number of distinct irreducible factors 
of the polynomials $\{P_n(t)\}_{n \in \N}$ is infinite. If $K$ has
characteristic $p>0$ and $q=p^e$ varies over the powers of $p$, then the
polynomials $\{P_{q-2}(t)\}_{q=p^e}$ have infinitely many distinct irreducible
factors.

Consequently the number of distinct irreducible factors of the homogeneous
polynomials $\{Q_n(s,t)\}_{n \in \N}$ as well as $\{Q_{q-2}(s,t)\}_{q=p^e}$
is also infinite.
\end{lemma}

\begin{proof}
It follows from Lemma~\ref{roots} that $\{P_n(t)\}_n$ as well as 
$\{P_{q-2}(t)\}_{q=p^e}$ have infinitely many distinct irreducible factors in
$K[t]$. 
\end{proof}

\section{Examples over integral domains}

We can now construct a domain which has a local cohomology module with
infinitely many associated primes:

\begin{theorem}\label{domain}
Let $K$ be an arbitrary field, and consider the integral domain 
$$
R=K[s,t,u,v,x,y]/\big(su^2x^2+tuxvy+sv^2y^2\big). 
$$
Then the local cohomology module $H^2_{(x,y)}(R)$ has infinitely many
associated prime ideals. Also, if we consider the local domain $R_{\m}$ where
$\m=(s,t,u,v,x,y)R$, then $H^2_{(x,y)}\left(R_{\m}\right)$ has infinitely many
associated primes.

If $\mathbb S$ is any infinite set of positive integers, then the set
$\bigcup_{m \in \mathbb S}\Ass R/(x^m,y^m)$ is infinite; in particular, if $K$
has characteristic $p>0$, then $\bigcup_{e \in \N}\Ass R/(x^{p^e},y^{p^e})$ is
infinite. The same conclusions hold if we replace the hypersurface $R$ by its
specialization $R/(s-1)$ or by the localization $R_{\m}$.
\end{theorem}

\begin{proof}
The assertions regarding local cohomology follow from Theorem~\ref{thmlc} and
Lemma~\ref{inf}. These, along with Proposition~\ref{contain}, imply the results
for generalized Frobenius powers of ideals; see also the remark below.
\end{proof}

\begin{remark} 
Specializing $s=1$ and working with the hypersurface
$$
S = R/(s-1) = K[t,u,v,x,y]/\big(u^2x^2+tuxvy+v^2y^2\big),
$$
similar arguments show that $H^2_{(x,y)}(S)$ has infinitely many associated
primes. This gives an example of a four dimensional integral domain $S$ for
which $H^2_{(x,y)}(S)$ has infinitely many associated prime ideals. However it
remains an open question whether a local cohomology module $H^i_{\a}(T)$ has
infinitely many associated primes where $T$ is a {\em local}\/ ring of dimension
four. This is of interest in view of Marley's results that the local cohomology
of a Noetherian local ring of dimension less than four has finitely many
associated primes, \cite{Ma}.
\end{remark}

For the assertion regarding the associated primes of generalized Frobenius
powers of an ideal, the hypersurface $R$ of Theorem~\ref{domain} can be modified
to obtain a three-dimensional local domain, or a two-dimensional non-local
domain:

\begin{theorem}\label{secondmain}
Let $K$ be an arbitrary field, and consider the integral domain
$$
A=K[s,t,x,y]/\big(sx^2+txy+sy^2\big). 
$$
Then the set $\bigcup_{n \in \N}\Ass A/(x^n,y^n)$ is infinite. The same
conclusion holds if we replace $A$ by the specialization $A/(s-1)$ or by the
localization $A_{(s,t,x,y)}$.
\end{theorem}

The proof of the theorem is again an immediate consequence of
Theorem~\ref{thmlc} and Lemma~\ref{inf}, but we feel it is of interest to
explicitly determine the infinite set $\bigcup_{n \in \N}\Ass A/(x^n,y^n)$ at
least in this one example, and we record the result as Theorem~\ref{genFrob}.
If $K$ has characteristic $p>0$, this theorem also shows that the set
$\bigcup_{e \in \N}\Ass A/(x^{p^e},y^{p^e})$ is finite. We next record some
preliminary computations which will be needed in determining the associated
primes of the ideals $(x^n,y^n)A$, and will also be used later in
\S\ref{freg-frat}.

\begin{lemma}\label{colon}
Consider the polynomial ring $K[s,t,x,y]$ and $m,n \ge 1$. Then
\begin{enumerate}
\item \quad $xy^{n-1}Q_{n-1} \ \in \ (x^n,y^n,sx^2+txy+sy^2)$,
\item \quad $(x^n,y^n,sx^2+txy+sy^2) : (s^mxy^{n-1}) \ = \ (x,y,Q_{n-1})$,
and
\item \quad $(x^n,y^n,x^2+txy+y^2) : (xy^{n-1}) \ = \ (x,y,P_{n-1})$.
\end{enumerate}
\end{lemma}

\begin{proof}
(1) The case $n = 1$ holds trivially. Using the equation 
$tQ_i=Q_{i+1}+s^2 Q_{i-1}$ for $1 \le i \le n-2$, we get
\begin{multline*}
(sx^2+txy+sy^2)(sx)^{n-2-i}y^iQ_i \ = \ s^{n-1-i}x^{n-i}y^iQ_i \ + \ 
s^{n-1-i}x^{n-2-i}y^{i+2}Q_i \\
+ \ s^{n-2-i}x^{n-1-i}y^{i+1}Q_{i+1} \ + \ s^{n-i}x^{n-1-i}y^{i+1}Q_{i-1},
\end{multline*}
and taking an alternating sum gives us
\begin{multline*}
\sum_{i=0}^{n-2} (-1)^i (sx^2+txy+sy^2)(sx)^{n-2-i}y^iQ_i \\
= \ s^{n-1}x^nQ_0 \ + \ (-1)^{n-2}xy^{n-1}Q_{n-1} \ + \ (-1)^{n-2}sy^nQ_{n-2}.
\end{multline*}
This shows that $xy^{n-1}Q_{n-1} \in (x^n,y^n,sx^2+txy+sy^2)$.

(2) If $n=1$ we have the unit ideal on each side of the asserted equality, so
we may assume $n \ge 2$ for the rest of this proof. It is easy to verify that 
$$
sxy^{n-1}(x,y) \ \subseteq \ (x^n,y^n,sx^2+txy+sy^2).
$$
Let $h \in K[s,t]$ be an element such that
$$
hs^mxy^{n-1} \ \in \ (x^n,y^n,sx^2+txy+sy^2).
$$
Using the grading where $\deg s=\deg t =0$ and $\deg x=\deg y=1$, there exist
elements $\alpha$, $\beta$, and $d_0,\dots,d_{n-2}$ in $K[s,t]$ with
\begin{multline*}
hs^m xy^{n-1} = (d_0x^{n-2}-d_1x^{n-3}y+\cdots+(-1)^{n-2}d_{n-2}y^{n-2}) 
(sx^2+txy+sy^2) \\
+ \alpha x^n+\beta y^n.
\end{multline*}
Comparing coefficients of $x^{n-1}y,x^{n-2}y^2,\dots,xy^{n-1}$, we get
\begin{align*}
s d_1 - t d_0 & = 0, \\ 
\qquad s d_{i+2} - t d_{i+1}+s d_i & = 0 \quad \text{for all} \quad
0 \le i \le n-4,\\
(-1)^{n-2} (t d_{n-2} - s d_{n-3}) &= hs^m.
\end{align*}
In particular,
$$
d_1 = (t/s)d_0 \quad \text{and} \quad d_{i+2} = (t/s)d_{i+1} - d_i \quad
\text{for all} \quad 0 \le i \le n-4,
$$ 
and consequently $d_i = d_0 P_i(t/s)$ for $0 \le i \le n-2$, where the $P_i$ 
are the polynomials defined recursively in \S\ref{d=2}. This gives us
$$
hs^m = (-1)^{n-2} \big(td_0P_{n-2}(t/s) - sd_0P_{n-3}(t/s)\big)
= (-1)^{n-2}sd_0P_{n-1}(t/s),
$$
and so $hs^{m+n-2} = (-1)^{n-2} d_0Q_{n-1}$. Since $s$ and $Q_{n-1}$ are 
relatively prime in $K[s,t]$, we see that $h$ is a multiple of $Q_{n-1}$. 

(3) is the inhomogeneous case of (2), and is left to the reader.
\end{proof}

\begin{lemma}\label{decomposition}
Let $A=K[s,t,x,y]/(sx^2+txy+sy^2)$, and $n \ge 1$ be an arbitrary integer.
\begin{enumerate}
\item For all $1 \le i \le n$, we have 
$s^{i-1}x^iy^{n-i} \in (x^n,y^n,xy^{n-1})$. In particular, 
$$
s^{n-1}(x,y)^n \subseteq (x^n,y^n,xy^{n-1}) \quad \text{and} \quad 
s^n(x,y)^n \subseteq (x^n,y^n,sxy^{n-1}).
$$
\item Also, $t^n(x,y)^n \subseteq (x^n,y^n,sxy^{n-1})$.
\item If $n \ge 2$, the ideal $(x^n,y^n,sxy^{n-1})$ has a primary decomposition
$$
(x^n,y^n,sxy^{n-1}) = (x,y)^n \cap (x^n,y^n,sxy^{n-1},s^n,t^n).
$$
\end{enumerate}
\end{lemma}

\begin{proof}
For (1) we use induction on $i$ to show that $s^{i-1}x^iy^{n-i} \in 
(x^n,y^n,xy^{n-1})$. This is certainly true if $i=1$ and, assuming the result 
for integers less than $i$, observe that
\begin{align*}
s^{i-1}x^iy^{n-i} &= -s^{i-2} x^{i-2} y^{n-i} (txy+sy^2)
= -s^{i-2} t x^{i-1} y^{n-i+1} - s^{i-1} x^{i-2} y^{n-i+2} \\
&\in (x^n,y^n,xy^{n-1}).
\end{align*}
Next, the equation $txy=-(sx^2+sy^2)$ gives us
$$
t(x,y)^n \subseteq (x^n,y^n)+s(x,y)^n,
$$
and using this inductively, we get
$$
t^n(x,y)^n \subseteq (x^n,y^n)+s^n(x,y)^n \subseteq (x^n,y^n,sxy^{n-1}),
$$
which proves (2). 

We next use the grading on the hypersurface $A$ where $\deg s=\deg t=0$ and 
$\deg x=\deg y=1$. If $\alpha$ and $\beta$ are nonzero homogeneous elements of
$A$ with $\alpha s^n+\beta t^n \in (x,y)^n$, then $\alpha$ and $\beta$ must
have degree at least $n$, and therefore belong to the ideal $(x,y)^n$. This
shows that
$$
(s^n,t^n) \cap (x,y)^n = (s^n,t^n)(x,y)^n,
$$
and using (1) and (2) we get
$$
(s^n,t^n) \cap (x,y)^n \subseteq (x^n,y^n,sxy^{n-1}).
$$

The intersection asserted in (3) follows immediately from this, and it remains
to verify that the ideals $\q_1=(x,y)^n$ and $\q_2=(x^n,y^n,sxy^{n-1},s^n,t^n)$
are indeed primary ideals. The radical of $\q_2$ is the maximal ideal
$(s,t,x,y)$, so $\q_2$ is a primary ideal. Using the earlier grading, any
homogeneous zerodivisor in the ring $A/\q_1$ must have positive degree, and
hence must be nilpotent. Consequently $\q_1$ is a primary ideal as well.
\end{proof}

\begin{theorem}\label{genFrob}
Let $A=K[s,t,x,y]/(sx^2+txy+sy^2)$ where $K$ is a field. Then
$\Ass A/(x^2,y^2)=\{(x,y),(t,x,y)\}$ and
$$
\Ass \frac{A}{(x^n,y^n)}=\{(x,y), (s,t,x,y)\} \cup
\Ass \frac{A}{(x,y,Q_{n-1})} \quad \text{for} \quad n \ge 3.
$$
In particular, $\bigcup_{n \in \N}\Ass A/(x^n,y^n)$ is an infinite set. If $K$
is an algebraically closed field, let
$$
\mathcal S =\left\{\left(x,y,t-s\xi-s\xi^{-1}\right)A \ | \ \xi\in K, \ \xi^n=1
\text{ for some } n \ge 1, \text{ and } \xi \neq \pm1 \right\}. 
$$
In the case that $K$ has characteristic zero,
$$
\bigcup_{n \ge 1}\Ass \frac{A}{(x^n,y^n)} =\{(x,y),(t,x,y),(s,t,x,y)\} \cup 
\mathcal S,
$$
and if $K$ has positive characteristic, then
$$
\bigcup_{n \ge 1}\Ass \frac{A}{(x^n,y^n)} 
= \{(x,y),(t,x,y),(s,t,x,y),(t-2s,x,y),(t+2s,x,y)\} \cup \mathcal S.
$$
\end{theorem}

\begin{proof}
It is easily checked that $(x,y)^2 \cap (x^2,y^2,t)$ is a primary decomposition
of $(x^2,y^2)$, so we need to compute $\Ass A/(x^n,y^n)$ for $n\ge 3$. By 
Lemma~\ref{colon}\,$(2)$ we have an exact sequence
$$
0 \xrightarrow{\phantom{\cdot sxy^{n-1}}} \frac{A}{(x,y,Q_{n-1})} 
\xrightarrow{\cdot sxy^{n-1}} \frac{A}{(x^n,y^n)} 
\xrightarrow{\phantom{\cdot sxy^{n-1}}} \frac{A}{(x^n,y^n,sxy^{n-1})} 
\xrightarrow{\phantom{\cdot sxy^{n-1}}} 0,
$$
and consequently
$$
\Ass \frac{A}{(x,y,Q_{n-1})} \subseteq\Ass \frac{A}{(x^n,y^n)} \subseteq 
\Ass \frac{A}{(x,y,Q_{n-1})} \cup\Ass \frac{A}{(x^n,y^n,sxy^{n-1})}.
$$
By Lemma~\ref{decomposition}
$\Ass A/(x^n,y^n,sxy^{n-1}) = \{(x,y), (s,t,x,y)\}$, and so it suffices to
verify that the prime ideals $\p_1=(x,y)$ and $\p_2=(s,t,x,y)$ are indeed 
associated primes of $A/(x^n,y^n)$. This follows since $\p_1$ is a minimal 
prime of $(x^n,y^n)$ and 
$$
\p_2 = (x^n,y^n):(xy)^{n-1}.
$$

If $K$ is an algebraically closed field, the polynomials $Q_i(s,t)$ split into
linear factors determined by the roots of $P_i(t)$, which are computed in
Lemma~\ref{roots}.
\end{proof}

\begin{remark}\label{examplelcfrob}
If $K$ is a field of characteristic $p>0$ and $A$ is the hypersurface
$$
A=K[s,t,x,y]/(sx^2+txy+sy^2)
$$
as above, we saw that the set $\bigcup_{n \in \N}\Ass A/(x^n,y^n)$ is infinite.
However, the set $\bigcup_{e \in \N}\Ass A/(x^{p^e},y^{p^e})$ is finite since,
by Theorem~\ref{genFrob},
$$
\Ass \frac{A}{(x^{p^e},y^{p^e})}=\{(x,y), (s,t,x,y)\} \cup 
\Ass \frac{A}{(x,y,Q_{p^e-1})} \quad \text{for} \quad p^e \ge 3,
$$
and $Q_{p^e-1}$ is a power of $Q_{p-1}$ by Proposition~\ref{threediag}. Using
Proposition~\ref{contain}, the set $\Ass H^2_{(x,y)}(R)$ is finite as well.
Consequently we have a strict inclusion
$$
\Ass H^2_{(x,y)}(R) \subsetneq \bigcup_{n \in \N} 
\Ass \frac{R}{(x^n,y^n)R}.
$$
The set $\bigcup_{n \in \N}\Ass A/(x^n,y^n)$ has been explicitly computed in
Theorem~\ref{genFrob}, and we next observe that the only associated prime of
$H^2_{(x,y)}(R)$ is the maximal ideal $\m=(s,t,x,y)$. The module
$H^2_{(x,y)}(R)$ is generated over $R$ by the elements $\eta_q = [1+(x^q,y^q)]$
for $q=p^e$, and it suffices to show that $\eta_q$ is killed by a power of
$\m$. It is immediately seen that $x^q$ and $y^q$ kill $\eta_q$, and for the
remaining cases note that
$$
s^q \eta_q = [s^q x^{2q} + (x^{3q},y^q)] = 0 \quad \text{and} \quad
t^q \eta_q = [t^q x^qy ^q + (x^{2q},y^{2q})] = 0.
$$
\end{remark}

\section{F-regular and unique factorization domain examples}\label{freg-frat}

In Theorem~\ref{domain} we proved that for the hypersurface
$$
R=K[s,t,u,v,x,y]/\big(su^2x^2+tuxvy+sv^2y^2\big), 
$$
the local cohomology module $H^2_{(x,y)}(R)$ has infinitely many associated
prime ideals. This ring $R$, while a domain, is not normal. In
Theorem~\ref{thmFratex} we construct examples over normal hypersurfaces, in
fact over hypersurfaces of characteristic zero with rational singularities, as
well as over F-regular hypersurfaces of positive characteristic. F-regularity
is a notion arising from the theory of tight closure developed by Hochster and
Huneke in \cite{HHjams}. A brief discussion may be found in \S\ref{tc}, though
for details of the theory and its applications, we refer the reader to
\cite{HHjams, HHbasec, HHjalg} and \cite{Hu2}.

\begin{theorem}\label{thmFratex}
Let $K$ be an arbitrary field, and consider the hypersurface
$$
S=\frac{K[s,t,u,v,w,x,y,z]}{\big(su^2x^2+sv^2y^2+tuxvy+tw^2z^2\big)}.
$$
Then $S$ is a normal domain for which the local cohomology module
$H^3_{(x,y,z)}(S)$ has infinitely many associated prime ideals. This is
preserved if we replace $S$ by $S/(s-1)$ or by the localization
$S_{(s,t,u,v,w,x,y,z)}$. If $K$ has characteristic zero, then $S$ has rational
singularities, and if $K$ has characteristic $p>0$, then $S$ is F-regular.
\end{theorem}

\begin{proof}
We defer the proof that $S$ has rational singularities or is F-regular, see
Lemma~\ref{freg} below. Normality follows from this, or may be proved directly
using the Jacobian criterion. Let $B$ be the subring of $S$ generated, as a
$K$-algebra, by the elements $s,t,a=ux,b=vy$ and $c=wz$, i.e.,
$$
B=K[s,t,a,b,c]/\big(sa^2+sb^2+tab+tc^2\big).
$$
For integers $n \ge 1$, let
$$
\eta_n = \left[s(ux)(vy)^{n-1}+(x^n,y^n,z)\right] \in H^3_{(x,y,z)}(S).
$$
Using $S_0=K[s,t]$ as the subring of $S$ of elements of degree zero,
Proposition~\ref{multigraded} implies that
$$
{\ann}_{S_0} \eta_n = (a^n,b^n,c)B :_{S_0} sab^{n-1}
$$
and then Lemma~\ref{colon}\,$(2)$ give us
$$
(a^n,b^n,c)B :_{S_0} sab^{n-1} = (Q_{n-1})S_0,
$$
where the $Q_i$ are the polynomials defined recursively in \S\ref{d=2}. Using
Lemma~\ref{subring} and Lemma~\ref{inf}, it follows that $H^3_{(x,y,z)}(S)$
has infinitely many associated prime ideals.
\end{proof}

It remains to prove that the hypersurface $S$ in Theorem~\ref{thmFratex} has
rational singularities or is F-regular, depending on the characteristic. The
results of \cite{SW} provide a direct proof that the hypersurface $S$ has
rational singularities in characteristic zero. However, instead of relying on
this, we prove here that if $K$ has positive characteristic, then $S$ is
F-regular. Using \cite[Theorem 4.3]{sm-ratsing}, it then follows that $S$ has
rational singularities when $K$ has characteristic zero. We first record an
elementary lemma:

\begin{lemma}\label{socle}
Let $(S,\m)$ be an $\N$-graded Gorenstein domain of dimension $d$, finitely
generated over a field $[S]_0=K$ of characteristic $p>0$, and let
$\eta \in H^d_{\m}(S)$ denote a socle generator. Let $c \in R$ be a nonzero
element such that $S_c$ is regular. Then $S$ is F-regular if and only if there
exists an integer $e \ge 1$ such that $\eta$ belongs to the $S$-span of
$cF^e(\eta)$.
\end{lemma}

\begin{proof}
If $S$ is F-regular then the zero submodule of $H^d_{\m}(S)$ is tightly closed,
i.e., $0^*_{H^d_\m(S)}=0$, and so there exists a positive integer $e$ such that
$cF^e(\eta) \neq 0$. Since $\eta$ generates the socle of $H^d_{\m}(S)$, which
is one-dimensional, $\eta$ must belong to the $S$-span of $cF^e(\eta)$.

Conversely, assume that $\eta$ belongs to the $S$-span of $cF^e(\eta)$ for some
$e \ge 1$. Then $cF^e(\eta) \neq 0$, and so the Frobenius morphism $F:
H^d_\m(S) \longrightarrow H^d_\m(S)$ is injective. It follows from
\cite[Proposition~6.11]{HRinv} that the ring $S$ is F-pure. By
\cite[Theorem~6.2]{HHbasec}, the element $c$ has a power which is a test
element but then, since $S$ is F-pure, $c$ itself must be a test element. The
condition $c F^e(\eta) \not = 0$ implies that $\eta \notin 0^*_{H^d_\m(S)}$.
Consequently $0^*_{H^d_\m(S)}=0$, and it follows that $S$ is F-regular.
\end{proof}

\begin{lemma}\label{freg}
Let $K$ be a field and consider the hypersurface
$$
S=\frac{K[s,t,u,v,w,x,y,z]}{\big(su^2x^2+sv^2y^2+tuxvy+tw^2z^2\big)}.
$$
If $K$ has characteristic $p>0$, then $S$ is F-regular. If $K$ has
characteristic zero, then $S$ has rational singularities.
\end{lemma}

\begin{proof}
We first consider the case where $K$ has characteristic $p>0$. It is easily
checked that $S_{twz}$ is a regular ring. We may compute $H^7_{\m}(S)$ using
the \v Cech complex with respect to the system of parameters
$s,u,x,v,y,w-t,z-t$. The socle of $H^7_{\m}(S)$ is spanned by the element
$$
\eta = [t^4 +(s,u,x,v,y,w-t,z-t)] \in H^7_{\m}(S).
$$
Since $S_{twz}$ is regular it suffices, by Lemma~\ref{socle}, to show that
$\eta$ belongs to the $S$-span of $twzF^e(\eta)$ for some $e\ge 1$, i.e., that
$$
t^4 \big(suxvy(w-t)(z-t) \big)^{q-1} \in 
\big(twzt^{4q},s^q,u^q,x^q,v^q,y^q,(w-t)^q,(z-t)^q \big)S \qquad (*) 
$$
for some $q=p^e$. We shall consider here the case $p \ge 5$, and the interested
reader may verify that $(*)$ holds with $q = 2^3$ and $q = 3^2$ in the
remaining cases $p=2$ and $p=3$ respectively. It suffices to show that
$$
t^4(suxvy)^{p-1} \in \big(t^{4p+3},s^p,u^p,x^p,v^p,y^p,w-t,z-t\big)S.
$$
Working in the polynomial ring $A=K[s,t,u,v,x,y]$, it is enough to check that
$t^4(suxvy)^{p-1} \in \a+(t^{5p-1})A$, where
$$
\a=(x^p,y^p,su^2x^2+sv^2y^2+tuxvy+t^5)A.
$$
We observe that
\begin{align*}
t^{5p-1} & \equiv t^4 \big(su^2x^2+sv^2y^2+tuxvy \big)^{p-1} \mod \a \\
& = t^4 \sum_{i,j} \binom{p-1}{i} \binom{p-1-i}{j} \big(su^2x^2 \big)^i
\big(sv^2y^2\big)^j \big(tuxvy\big)^{p-1-i-j} \mod \a \\
& = t^4 \sum_{i,j} \binom{p-1}{i}\binom{p-1-i}{j} s^{i+j} t^{p-1-i-j}
(ux)^{p-1+i-j} (vy)^{p-1-i+j} \mod \a.
\end{align*}
The only terms which contribute $\mod (x^p,y^p)$ are those for which $i=j$, and
so
$$
t^{5p-1} \equiv t^4 \sum_{i=0}^{(p-1)/2} \binom{p-1}{i}\binom{p-1-i}{i} 
s^{2i} t^{p-1-2i} (uxvy)^{p-1} \mod \a.
$$
When $2i < p-1$, the corresponding summand in the above expression is a
multiple of $t^5 (uxvy)^{p-1}$, which is an element of $\a$. Thus
$$
t^{5p-1} \equiv t^4 \binom{p-1}{(p-1)/2} s^{p-1} (uxvy)^{p-1} \mod \a.
$$
Since the binomial coefficient occurring above is a unit,
$t^4(suxvy)^{p-1} \in \a+(t^{5p-1})A$, which completes the proof that $S$ is
F-regular.

It remains to show that $S$ has rational singularities in the case $K$ has
characteristic zero. By \cite[Theorem~4.3]{sm-ratsing}, it suffices to show
that $S$ has F-rational type, i.e., that for all but finitely many prime
integers $p$, the fiber over $p\Z$ of the map
$$
\Z \longrightarrow 
\frac{\Z[s,t,u,v,w,x,y,z]}{\big(su^2x^2+sv^2y^2+tuxvy+tw^2z^2\big)}
$$
is an F-rational ring. While this is indeed true for {\it all}\/ prime integers
$p$, our earlier computation for $p \ge 5$ certainly suffices.
\end{proof}

We next construct unique factorization domains with similar behaviour.

\begin{theorem}
\label{lcufd}
Let $K$ be an arbitrary field, and consider the hypersurface
$$
T=\frac{K[r,s,t,u,v,w,x,y,z]}{\big(su^2x^2+sv^2y^2+tuxvy+rw^2z^2\big)}
$$
Then $T$ is a unique factorization domain for which the local cohomology module
$H^3_{(x,y,z)}(T)$ has infinitely many associated prime ideals. This is
preserved if we replace $T$ by the localization at its homogeneous maximal
ideal. The hypersurface $T$ has rational singularities if $K$ has
characteristic zero, and is F-regular in the case of positive characteristic.
\end{theorem}

\begin{proof}
It is easily verified that $T$ is a normal domain, in particular, the element
$t-r \in T$ is a nonzerodivisor. Note that
$$
T/(t-r) \cong K[s,t,u,v,w,x,y,z]/\big(su^2x^2+sv^2y^2+tuxvy+tw^2z^2\big)
$$
is F-regular or F-rational by Lemma~\ref{freg}. The rational singularity
property deforms by \cite{El}, and F-regularity deforms for Gorenstein rings by
\cite[Corollary~4.7]{HHbasec}. It follows that $T$ has rational singularities
if $K$ has characteristic zero, and is F-regular in positive characteristic.

We next prove that $T$ is a unique factorization domain. Consider the
multiplicative system $W \subset T$ generated by the elements $w$ and $z$.
Since $W$ is generated by prime elements, by Nagata's Theorem it suffices to
verify that $W^{-1}T$ is a unique factorization domain, see
\cite[Theorem~6.3]{Sa} or \cite[Corollary~7.3]{Fo}. But 
$$
W^{-1}T=K[s,t,u,v,x,y,w,w^{-1},z,z^{-1}]
$$
is a localization of a polynomial ring, and hence is a unique factorization
domain.

For integers $n \ge 1$, consider
$$
\eta_n = \left[s(ux)(vy)^{n-1}+(x^n,y^n,z)\right] \in H^3_{(x,y,z)}(T).
$$
As in the proof of Theorem~\ref{thmFratex}, we use Proposition~\ref{multigraded}
and Lemma~\ref{colon}\,$(2)$ to compute ${\ann}_{T_0}\eta_n$ where
$T_0 = K[r,s,t]$. Setting $a=ux$, $b=vy$, and $c=wz$, we see that
$$
{\ann}_{T_0} \eta_n = (Q_{n-1})T_0.
$$
By Lemma~\ref{subring} and Lemma~\ref{inf}, it follows that $H^3_{(x,y,z)}(T)$
has infinitely many associated prime ideals.
\end{proof}

\section{An application to tight closure theory} \label{tc}

Let $R$ be a ring of characteristic $p>0$, and $R^\circ$ denote the complement
of the minimal primes of $R$. For an ideal $\a=(x_1,\dots,x_n)$ of $R$ and a
prime power $q=p^e$, we use the notation $\a^{[q]}=(x_1^q,\dots,x_n^q)$. The
{\it tight closure}\/ of $\a$ is the ideal
$$
\a^* = \{z \in R \ | \ \text{there exists } c\in R^\circ \text{ for which }
cz^q \in \a^{[q]} \text{ for all } q\gg 0 \},
$$
see \cite{HHjams}. A ring $R$ is {\it F-regular}\/ if $\a^* = \a$ for all
ideals $\a$ of $R$ and its localizations. 

More generally, let $F$ denote the Frobenius functor, and $F^e$ its $e$\/th
iteration. If an $R$-module $M$ has presentation matrix $(a_{ij})$, then
$F^e(M)$ has presentation matrix $(a_{ij}^q)$, where $q=p^e$. For modules
$N \subseteq M$, we use $N^{[q]}_M$ to denote the image of $F^e(N) \to F^e(M)$. 
We say that an element $m \in M$ is in the {\it tight closure of $N$ in $M$},
denoted $N^*_M$, if there exists an element $c \in R^\circ$ such that
$cF^e(m) \in N^{[q]}_M$ for all $q \gg 0$. While the theory has found several
applications, the question whether tight closure commutes with localization
remains open even for finitely generated algebras over fields of positive
characteristic.

Let $W$ be a multiplicative system in $R$, and $N \subseteq M$ be finitely
generated $R$-modules. Then
$$
W^{-1}(N^*_M) \subseteq (W^{-1}N)^*_{W^{-1}M},
$$
where $W^{-1}(N^*_M)$ is identified with its image in $W^{-1}M$. When this
inclusion is an equality, we say that {\it tight closure commutes with
localization at $W$ for the pair $N \subseteq M$}. It may be checked that
this occurs if and only if tight closure commutes with localization at $W$
for the pair $0 \subseteq M/N$. Following \cite{AHH}, we set
$$
G^e(M/N) = F^e(M/N)/(0^*_{F^e(M/N)}).
$$
An element $c \in R^\circ$ is a {\it weak test element}\/ if there exists
$q_0=p^{e_0}$ such that for every pair of finitely generated modules
$N \subseteq M$, an element $m \in M$ is in $N^*_M$ if and only if
$cF^e(m) \in N^{[q]}_M$ for all $q \ge q_0$. By \cite[Theorem 6.1]{HHbasec}, if
$R$ is of finite type over an excellent local ring, then $R$ has a weak test
element. 

\begin{prop}\cite[Lemma 3.5]{AHH}
Let $R$ be a ring of characteristic $p>0$ and $N \subseteq M$ be finitely
generated $R$-modules. Then the tight closure of $N \subseteq M$ commutes with
localization at any multiplicative system $W$ which is disjoint from the set
$\bigcup_{e\in \N}\Ass F^e(M)/N^{[q]}_M$.

If $R$ has a weak test element, then the tight closure of $N \subseteq M$ also
commutes with localization at multiplicative systems $W$ disjoint from the set
$\bigcup_{e\in \N}\Ass G^e(M/N)$.
\end{prop}

Consider a bounded complex $\Pdot$ of finitely generated projective $R$-modules,
$$
0 \longrightarrow P_n \overset{d_n}\longrightarrow P_{n-1} \longrightarrow
\cdots \overset{d_1}\longrightarrow P_0 \longrightarrow 0.
$$
The complex $\Pdot$ is said to have {\it phantom homology}\/ at the $i$\/th
spot if
$$
\Ker d_i \subseteq (\Im d_{i+1})^*_{P_i}.
$$
The complex $\Pdot$ is {\it stably phantom acyclic}\/ if $F^e(\Pdot)$ has
phantom homology at the $i$\/th spot for all $i\ge 1$, for all $e \ge 1$. An
$R$-module $M$ has {\it finite phantom projective dimension}\/ if there exists
a bounded stably phantom acyclic complex $\Pdot$ of projective $R$-modules,
with $H_0(\Pdot) \cong M$.

\begin{theorem}\cite[Theorem 8.1]{AHH}
Let $R$ be an equidimensional ring of positive characteristic, which is of
finite type over an excellent local ring. If $N \subseteq M$ are finitely
generated $R$-modules such that $M/N$ has finite phantom projective dimension,
then the tight closure of $N$ in $M$ commutes with localization at $W$ for
every multiplicative system $W$ of $R$.
\end{theorem}

The key points of the proof are that for $M/N$ of finite phantom projective
dimension, the set $\bigcup_e\Ass G^e(M/N)$ has finitely many maximal elements,
and that if $(R,\m)$ is a local ring, then there a positive integer $B$ such
that for all $q=p^e$, the ideal $\m^{Bq}$ kills the local cohomology module
$$
H^0_\m \left(G^e(M/N) \right).
$$
For more details on this approach to the localization problem, we refer the
reader to the papers \cite{AHH, Ho, Ka1, SN}, and \cite[\S12]{Hu2}.
Specializing to the case where $M=R$ and $N=\a$ is an ideal, we note that
$$
G^e(R/\a) \cong R/ \big(\a^{[q]}\big)^*, \quad \text{where} \quad q=p^e.
$$
This raises the questions:

\begin{question}\cite[page 90]{Ho}\label{frob}
Let $R$ be a Noetherian ring of characteristic $p>0$, and $\a$ an ideal of $R$.
\begin{enumerate}
\item Does the set $\bigcup_q\Ass R/ \a^{[q]}$ have finitely many maximal
elements?
\item Does $\bigcup_q\Ass R/ \big(\a^{[q]}\big)^*$ have finitely many
maximal elements?
\item For a complete local domain $(R,\m)$ and an ideal $\a \subset R$, is there
a positive integer $B$ such that
$$
\m^{Bq} H^0_\m \left(R/ (\a^{[q]})^* \right) = 0 \quad \text{for all} \quad
q=p^e \,?
$$
\end{enumerate}
\end{question}

Katzman proved that affirmative answers to Questions~\ref{frob}\,$(2)$ and
\ref{frob}\,$(3)$ imply that tight closure commutes with localization:

\begin{theorem}\cite{Ka1}
Assume that for every local ring $(R,\m)$ of characteristic $p>0$ and ideal
$\a \subset R$, the set $\,\bigcup_q\Ass R/\big(\a^{[q]}\big)^*$ has finitely
many maximal elements. Also, if for every ideal $\a \subset R$, there exists a
positive integer $B$ such that $\m^{Bq}$ kills 
$$
H^0_\m \left(R/ (\a^{[q]})^* \right) \quad \text{for all} \quad q=p^e,
$$
then tight closure commutes with localization for all ideals in Noetherian
rings of characteristic $p>0$.
\end{theorem}

These issues are, of course, the source of our interest in associated primes of
Frobenius powers of ideals. It should be mentioned that the situation for
{\it ordinary}\/ powers is well-understood: $\bigcup_n\Ass R/\a^n$ is finite
for any Noetherian ring $R$, see \cite{Br} or \cite{Ra}. However, for
{\it Frobenius} powers, Katzman showed that the maximal elements of
$\bigcup_q\Ass R/ \a^{[q]}$ need not form a finite set, thereby settling
Question~\ref{frob}\,$(1)$. We recall the example from \cite{Ka1}, discussed
earlier in Remark~\ref{moty}: if
$$
A = K[t,x,y]/\big(xy(x-y)(x-ty)\big),
$$
then the set $\bigcup_q\Ass R/ (x^q,y^q)$ is infinite. In this example
$(x^q,y^q)^*=(x,y)^q$ for all $q=p^e$ and so, in contrast,
$\bigcup_q\Ass A/(x^q,y^q)^*$ is finite.

\begin{remark}
In Theorem~\ref{thmFratex} we constructed an F-regular ring $S$ for which the
set $\Ass H^3_{(x,y,z)}(S)$ is infinite. By Proposition~\ref{contain} we have
$$
\Ass H^3_{(x,y,z)}(S) \subseteq \bigcup_{q=p^e}\Ass S/(x^q,y^q,z^q),
$$
and it follows that $\bigcup_{q}\Ass S/(x^q,y^q,z^q)$ must be infinite. Since
$S$ is F-regular, we have $(x^q,y^q,z^q)^*=(x^q,y^q,z^q)$ for all $q=p^e$, and
so $\bigcup_{q}\Ass S/(x^q,y^q,z^q)^*$ is infinite. The question remains
whether $\bigcup_q\Ass R/ \big(\a^{[q]}\big)^*$ has finitely many
{\it maximal}\/ elements for arbitrary rings $R$ of characteristic $p>0$, and
we next show that this has a negative answer as well, thereby settling
Question~\ref{frob}\,$(2)$. \end{remark}

\begin{theorem}\label{assmax}
Let $K$ be a field of characteristic $p>0$, and consider
$$
R=\frac{K[t,u,v,w,x,y,z]}{\big(u^2x^2+v^2y^2+tuxvy+tw^2z^2\big)}.
$$
Then $R$ is an F-regular ring, and the set
$$
\bigcup_{e \in \N}\Ass \frac{R}{\big(x^{p^e},y^{p^e},z^{p^e}\big)}
=\bigcup_{e \in \N}\Ass \frac{R}{\big(x^{p^e},y^{p^e},z^{p^e}\big)^*}
$$
has infinitely many maximal elements.
\end{theorem}

\begin{proof}
By Lemma~\ref{freg} the hypersurface
$$
S= K[s,t,u,v,w,x,y,z]/\big(su^2x^2+sv^2y^2+tuxvy+tw^2z^2\big)
$$
is F-regular, and therefore so is its localization
$$
S_s=\frac{K[t/s,u,v,w,x,y,z,s,1/s]}
{\big(u^2x^2+v^2y^2+(t/s)uxvy+(t/s)w^2z^2\big)}.
$$
The ring $S_s$ has a $\Z$-grading where $\deg s=1$, $\deg 1/s=-1$, and
the remaining generators, $t/s,u,v,w,x,y,z$, have degree $0$. By
\cite[Proposition~4.12]{HHjams} a direct summand of an F-regular ring is
F-regular, and so $R \cong {[S_s]}_0$ is F-regular.

For $q=p^e$, consider the ideals of $R$,
$$
\a_q \ = \ (x^q,y^q,z^q)R :_R t^quv^{q-2}x^2y^{q-1}z^{q-1}.
$$
Let $R_0=K[t]$. As in the proof of Theorem~\ref{thmFratex}, we may use
Proposition~\ref{multigraded} and Lemma~\ref{colon}\,$(3)$ to verify that
\begin{align*}
\a_q \cap R_0
&\ = \ (x^q,y^q,z^q)R :_{R_0} t^q(ux)(vy)^{q-2} xyz^{q-1} \\
&\ = \ (x^{q-1},y^{q-1},z)R :_{R_0} t^q(ux)(vy)^{q-2}
\ = \ P_{q-2} :_{R_0} t^q,\\
\end{align*}
where the $P_i$ are the polynomials defined recursively in \S\ref{d=2}. In
particular, this shows that $\a_q \neq R$ for $q \gg 0$. It is immediately
seen that $x,y,z \in \sqrt{\a_q}$, and we claim that $u,v,w \in \sqrt{\a_q}$.
To see that $u \in \a_q$, note that
$$
u(t^quv^{q-2}x^2y^{q-1}z^{q-1}) \ = \
t^q(u^2x^2)v^{q-2}y^{q-1}z^{q-1} \in (y^q,z^q).
$$
Next, observe that
$$
(vy)^2 \in ux(ux,vy)R+zR, \quad \text{and so} \quad 
(vy)^{q-1} \in (ux)^{q-2}(ux,vy)R+zR. 
$$
Using this,
$$
v(t^quv^{q-2}x^2y^{q-1}z^{q-1}) \ = \ t^q(vy)^{q-1}ux^2z^{q-1} \in (x^q,z^q),
$$
and so $v \in \a_q$. Finally, it is easily verified that $w^{q-1} \in \a_q$,
i.e., that
$$
w^{q-1}(st^quv^{q-2}x^2y^{q-1}z^{q-1}) \in (x^q,y^q,z^q),
$$
since $t^q(wz)^{q-1} \ \in \ (x^{q-2},y)$. We have now established
$$
{\Min}(\a_q) \ = \ {\Min}\big((u,v,w,x,y,z)R+(P_{q-2} :_{R_0} t^q)R \big),
$$
and so the minimal primes of $\a_q$ are maximal ideals of $R$. By
Lemma~\ref{inf} the union $\bigcup_q{\Min}(\a_q)$ is an infinite set, and so we
conclude that $\bigcup_q\Ass R/\big(x^q,y^q,z^q\big)$ has infinitely many
maximal elements.
\end{proof}

\begin{remark}
We would like to point out that the ring $R$ in Theorem~\ref{assmax} is a unique
factorization domain if $K = {\Z}/p{\Z}$ where $p$ is a prime with
$p \equiv 3 \mod 4$ or, more generally, if $K$ does not contain a square-root of
$-1$. In this case the polynomial $u^2x^2 + v^2y^2$ is irreducible, so
$f=uxvy+w^2z^2 \in R$ is a prime element. The ring $R_f$ is a localization of
$K[u,v,w,x,y,z]$, hence is a unique factorization domain. By Nagata's Theorem,
it follows that $R$ is a unique factorization domain.

For examples which do not depend on the field $K$, the interested reader may
verify that
$$
S=K(r)[t,u,v,w,x,y,z]/\big(u^2x^2+v^2y^2+tuxvy+rw^2z^2\big)
$$
is an F-regular unique factorization domain for which the set
$$
\bigcup_{e \in \N}\Ass \frac{S}{\big(x^{p^e},y^{p^e},z^{p^e}\big)}
=\bigcup_{e \in \N}\Ass \frac{S}{\big(x^{p^e},y^{p^e},z^{p^e}\big)^*}
$$
has infinitely many maximal elements.
\end{remark}

\section{Examples of small dimension}\label{d=4}

We analyze multidiagonal matrices with $d=4$ and use these computations to
obtain low-dimensional examples of integral domains of characteristic $p>0$
where the set of associated primes of Frobenius powers of an ideal is infinite.
The example in Theorem~\ref{domain}, after specializing $s=1$, is an integral
domain of dimension four. We construct here a hypersurface $A$ of dimension
two, which is an integral domain, and has an ideal $(x, y)A$ for which
$\bigcup_e\Ass A/\big(x^{p^e},y^{p^e}\big)$ is infinite. In view of
Proposition~\ref{threediag}, to construct such an example using
Theorem~\ref{thmlc}, we need to consider multidiagonal matrices with $d \ge 4$.

We start with the polynomial ring $A_0 = K[t]$ over a field $K$. Let $d = 4$,
and consider the matrices $M_n$ of multidiagonal form with respect to
$r_0=r_4=1, r_2=t$, and $r_1=r_3=0$, i.e.,
$$
M_n=\left[\begin{matrix}
t & 0 & 1 \cr
0 & t & 0 & 1 \cr
1 & 0 & t & 0 & 1 \cr
& \ddots & \ddots & \ddots & \ddots & \ddots \cr
& & 1 & 0 & t & 0 & 1 \cr
& & & 1 & 0 & t & 0 \cr
& & & & 1 & 0 & t \cr
\end{matrix}\right].
$$

We again use the convention $\det M_0=1$, and have $\det M_1=t$, $\det M_2=t^2$,
$\det M_3=t^3-t$, and the recursion
$$
\det M_{n+4} = t \det M_{n+3} - t \det M_{n+1}+\det M_n \quad \text{for all}
\quad n \ge 0.
$$
Using this, the generating function for $\det M_n$ is easily computed to be
$$
G(x) = \sum_{n \ge 0} \left(\det M_n\right) x^n = \frac{1}{1-tx+tx^3-x^4}
= \frac{1}{(1-x)(1+x)(1-tx+x^2)}.
$$
Set $F_n(t)=\det M_n$, which is a monic polynomial of degree $n$. We need to
analyze the distinct irreducible factors of the polynomials $\{F_n(t)\}$.

\begin{lemma}\label{roots2}
Let $K$ be an algebraically closed field, and consider the polynomials
$F_n(t)=\det M_n \in K[t]$ as above.
\begin{enumerate}
\item Let $\xi$ be a nonzero element of $K$ with $\xi \neq \pm 1$. If $n$ is an
odd integer, then
$$
F_n(\xi+\xi^{-1}) = \frac{(\xi^{n+3}-1)(\xi^{n+1}-1)}{\xi^n(\xi^2-1)^2},
$$
and so $F_n(\xi+\xi^{-1}) = 0$ if and only if $\xi^{n+3}=1$ or $\xi^{n+1}=1$.

\item If $n$ is an odd integer and $(n+3)(n+1)$ is invertible in $K$, then the
polynomial $F_n(t)$ has $n$ distinct roots of the form $\xi+\xi^{-1}$, where
$\xi \neq \pm 1$, and either $\xi^{n+3}=1$ or $\xi^{n+1}=1$.

\item If the characteristic of $K$ is an odd prime $p$, then $F_{q-2}(t)$ has 
$q-2$ distinct roots for all $q=p^e$.
\end{enumerate}
\end{lemma}

\begin{proof}
(1) Consider the generating function for the polynomials $F_n(t)$, 
$$
G(x)=\sum_{n \ge 0} F_n(t)x^n = \frac{1}{(1-x)(1+x)(1-tx+x^2)} \in K[t][[x]].
$$
If $\xi \in K$ with $\xi\neq 0$ and $\xi \neq \pm 1$, then
\begin{multline*}
\sum_{n \ge 0} F_n(\xi+\xi^{-1})x^n 
= \frac{1}{(1-x)(1+x)(1- \xi x)(1- \xi^{-1} x)} \\
= \frac{\sum x^n}{2(2-\xi-\xi^{-1})} 
+ \frac{\sum (-x)^n}{2(2+\xi+\xi^{-1})}
+ \frac{\xi^3 \sum (\xi x)^n}{(\xi^2-1)(\xi-\xi^{-1})} 
+ \frac{\xi^{-3} \sum (\xi^{-1} x)^n}{(\xi^{-2}-1)(\xi^{-1}-\xi)}.
\end{multline*}
Comparing the coefficients of $x^n$ and simplifying, we obtain the asserted
formula for $F_n(\xi+\xi^{-1})$.

(2) As we observed earlier in the proof of Lemma~\ref{roots}\,$(2)$, $\xi
+\xi^{-1} = \eta+\eta^{-1}$ if and only if $\xi$ equals $\eta$ or $\eta^{-1}$.
The only common roots of the polynomials $X^{n+3}-1=0$ and $X^{n+1}-1=0$ are
$\pm 1$. Since $n+3$ is invertible in the field $K$, the polynomial
$X^{n+3}-1=0$ has $n+1$ distinct roots $\xi$ with $\xi \neq \pm 1$. These give
the $(n+1)/2$ distinct roots $\xi+\xi^{-1}$ of $F_{n}(t)$. Similarly, the roots
of $X^{n+1}-1=0$ contribute $(n-1)/2$ other distinct roots of $F_{n}(t)$. But
then we have $(n+1)/2+(n-1)/2 = n$ distinct roots of the degree $n$ polynomial
$F_n(t)$ which, then, must be all its roots.

(3) Since $n=q-2$ is odd and $(n+3)(n+1)=(q+1)(q-1)$ is invertible in $K$, it
follows from $(2)$ that $F_{q-2}(t)$ has $q-2$ distinct roots.
\end{proof}

As a consequence of Lemma~\ref{roots2}, we immediately have:

\begin{lemma}\label{inf2}
Let $K$ be an arbitrary field of characteristic $p>2$. Then the polynomials
$\{F_{q-2}(t)\}_{q=p^e}$ have infinitely many distinct irreducible factors.
\end{lemma}

\begin{theorem}\label{lowdim}
Let $K$ be an arbitrary field of characteristic $p>2$, and consider the
integral domain 
$$
A=K[t,x,y]/\big(x^4+tx^2y^2+y^4\big). 
$$
Then the set $\bigcup_{e \in \N}\Ass A/\big(x^{p^e},y^{p^e}\big)$ is infinite.
\end{theorem}

\begin{proof}
The hypersurface $A$ arises from Theorem~\ref{thmlc} using the matrices $M_n$
of multidiagonal form with respect to $r_0=r_4=1, r_2=t$, and $r_1=r_3=0$. By
Lemma~\ref{inf2}, the set $\bigcup_e{\Min}\left(\det M_{p^e-2}\right)$ is
infinite, and so the result follows.
\end{proof} 

\section*{Acknowledgments} 

The first author is grateful to Mordechai Katzman and Uli Walther for
discussions regarding this material. The second author thanks Kamran
Divaani-Aazar and the Institute for Studies in Theoretical Physics and
Mathematics (IPM) in Tehran, Iran, for their interest and hospitality. Both
authors are grateful to Rodney Sharp for a careful reading of the manuscript,
and to The Mathematical Sciences Research Institute (MSRI) where parts of this
manuscript were prepared.


\begin{thebibliography}{BKS}

\bibitem[AHH]{AHH} I. M. Aberbach, M. Hochster, and C. Huneke, {\em
Localization of tight closure and modules of finite phantom projective
dimension}, J. Reine Angew. Math. {\bf 434} (1993), 67--114.

\bibitem[Br]{Br} M. Brodmann, {\em Asymptotic stability of ${\rm Ass} \,
(M/I\sp{n}M)$}, Proc. Amer. Math. Soc. {\bf 74} (1979), 16--18.

\bibitem[BF]{BF} M. P. Brodmann and A. Lashgari Faghani, {\em A finiteness
result for associated primes of local cohomology modules}, Proc. Amer. Math.
Soc. {\bf 128} (2000), 2851--2853. 

\bibitem[BKS]{BKS} M. P. Brodmann, M. Katzman, and R. Y. Sharp, {\em Associated
primes of graded components of local cohomology modules}, Trans. Amer. Math.
Soc. {\bf 354} (2002), 4261--4283.

\bibitem[BRS]{BRS} M. Brodmann, Ch. Rotthaus, and R. Y. Sharp, {\em On
annihilators and associated primes of local cohomology modules}, J. Pure Appl.
Alg. {\bf 153} (2000), 197--227.

\bibitem[El]{El} R. Elkik, {\em Singularit\'es rationnelles et d\'eformations},
Invent. Math. {\bf 47} (1978), 139--147.

\bibitem[Fo]{Fo} R. M. Fossum, {\em The divisor class group of a Krull domain},
Ergebnisse der Mathematik und ihrer Grenzgebiete {\bf 74}, Springer-Verlag,
New York-Heidelberg, 1973. 

\bibitem[He]{He} M. Hellus, {\em On the set of associated primes of a local
cohomology module}, J. Algebra {\bf 237} (2001), 406--419.

\bibitem[Ho]{Ho} M. Hochster, {\em The localization question for tight
closure}, in: Commutative algebra (International Conference, Vechta, 1994),
pp.~89--93, Vechtaer Universit\"atsschriften {\bf 13}, Verlag Druckerei Rucke
GmbH, Cloppenburg, 1994.

\bibitem[HH1]{HHjams} M. Hochster and C. Huneke, {\em Tight closure, invariant
theory, and the Brian\c con-Skoda theorem}, J. Amer. Math. Soc. {\bf 3} (1990),
31--116.

\bibitem[HH2]{HHbasec} M. Hochster and C. Huneke, {\em F-regularity, test
elements, and smooth base change}, Trans. Amer. Math. Soc. {\bf 346} (1994),
1--62.

\bibitem[HH3]{HHjalg} M. Hochster and C. Huneke, {\em Tight closure of
parameter ideals and splitting in module-finite extensions}, J. Algebraic Geom.
{\bf 3} (1994), 599--670. 

\bibitem[HH4]{HHchar0} M. Hochster and C. Huneke, {\em Tight closure in equal
characteristic zero}, in preparation.

\bibitem[HR]{HRinv} M. Hochster and J. L. Roberts, {\em Rings of invariants of
reductive groups acting on regular rings are Cohen-Macaulay}, Advances in Math.
{\bf 13} (1974), 115--175.

\bibitem[Hu1]{Hu1} C. Huneke, {\em Problems on local cohomology}, in: Free
resolutions in commutative algebra and algebraic geometry (Sundance, Utah,
1990), pp. 93--108, Res. Notes Math. {\bf 2}, Jones and Bartlett, Boston, MA,
1992.

\bibitem[Hu2]{Hu2} C. Huneke, {\em Tight closure and its applications}, CBMS
Regional Conference Series in Mathematics {\bf 88}, American Mathematical
Society, Providence, RI, 1996.

\bibitem[HS]{HS} C. L. Huneke and R. Y. Sharp, {\em Bass numbers of local
cohomology modules}, Trans. Amer. Math. Soc. {\bf 339} (1993), 765--779. 

\bibitem[Ka1]{Ka1} M.\ Katzman, {\em Finiteness of $\bigcup_e {\rm
Ass}\,F^e(M)$ and its connections to tight closure}, Illinois J. Math. {\bf 40}
(1996), 330--337.

\bibitem[Ka2]{Ka2} M. Katzman, {\em An example of an infinite set of associated
primes of a local cohomology module}, J. Algebra {\bf 252} (2002), 161--166. 

\bibitem[KS]{KS} K. Khashyarmanesh and Sh. Salarian, {\em On the associated
primes of local cohomology modules}, Comm. Alg. {\bf 27} (1999), 6191--6198.

\bibitem[Ly1]{Ly1} G. Lyubeznik, {\em Finiteness properties of local cohomology
modules (an application of $D$-modules to commutative algebra)}, Invent. Math.
{\bf 113} (1993), 41--55.

\bibitem[Ly2]{Ly2} G. Lyubeznik, {\em Finiteness properties of local cohomology
modules for regular local rings of mixed characteristic: the unramified case},
Comm. Alg. {\bf 28} (2000), 5867--5882.

\bibitem[Ly3]{Ly3} G. Lyubeznik, {\em Finiteness properties of local cohomology 
modules: a characteristic-free approach}, J. Pure Appl. Alg. {\bf 151} (2000), 
43--50. 

\bibitem[Ma]{Ma} T. Marley, {\em The associated primes of local cohomology
modules over rings of small dimension}, Manuscripta Math. {\bf 104} (2001),
519--525.

\bibitem[MV]{MV} T. Marley and J. C. Vassilev, {\em Cofiniteness and associated
primes of local cohomology modules}, J. Algebra {\bf 256} (2002), 180--193.

\bibitem[Ra]{Ra} L. J. Ratliff Jr., {\em On prime divisors of $I\sp{n}$, $n$
large}, Michigan Math. J. {\bf 23} (1976), 337--352.

\bibitem[Sa]{Sa} P. Samuel, {\em Lectures on unique factorization domains},
Tata Institute of Fundamental Research Lectures on Mathematics {\bf 30},
Bombay, 1964.

\bibitem[SN]{SN} R. Y. Sharp and N. Nossem, {\em Ideals in a perfect closure,
linear growth of primary decompositions, and tight closure}, Trans. Amer. Math.
Soc., posted on January 13, 2004, PII~S~0002-9947(04)03420-8 (to appear in
print).

\bibitem[Si]{Si} A. K. Singh, {\em $p$-torsion elements in local cohomology 
modules}, Math. Res. Lett. {\bf 7} (2000), 165--176.

\bibitem[SW]{SW} A. K. Singh and K.-i. Watanabe, {\em Multigraded rings,
rational singularities, and diagonal subalgebras}, in preparation.

\bibitem[Sm]{sm-ratsing} K. E. Smith, {\em F-rational rings have rational
singularities}, Amer. J. Math. {\bf 119} (1997), 159--180.

\bibitem[TZ]{TZ} R. Tajarod and H. Zakeri, {\em On the local-global principle
and the finiteness of associated primes of local cohomology modules}, Math. J.
Toyama Univ. {\bf 23} (2000), 29--40.

\end{thebibliography}
\end{document}